\documentclass[proc]{edpsmath}
\usepackage{path, graphicx}
\newtheorem{dfntnprpstn}[thrm]{Definition - Proposition}
\newtheorem{prprt}[thrm]{Property}

\newcommand{\A}[1]{{\vphantom{A}}_{#1}{A}}
\newcommand{\At}[1]{{\vphantom{A}}_{#1}{\widetilde{A}}}
\newcommand{\Id}{\mbox{Id}}
\newcommand{\nablat}{\widetilde{\nabla}}

\newcommand{\vectq}{\boldsymbol{q}}                                            
\newcommand{\vectf}{\boldsymbol{f}}
\newcommand{\vectv}{\boldsymbol{v}}

\newcommand{\vectWt}{\widetilde{\boldsymbol{W}}}
\newcommand{\vectW}{\boldsymbol{W}}
\newcommand{\vectY}{\boldsymbol{Y}}
\newcommand{\vectm}{\boldsymbol{m}}

\newcommand{\vecte}{\boldsymbol{e}}
\newcommand{\vectet}{\widetilde{\boldsymbol{e}}}
\newcommand{\vectx}{\boldsymbol{x}}
\newcommand{\vectxt}{\widetilde{\boldsymbol{x}}}

\newcommand{\trans}[1]{{#1}^{\texttt{t}}}
\newcommand{\contr}{{\;\boldsymbol :\;}}

\newcommand{\qy}{{q_y}}							
\newcommand{\qx}{{q_x}}
\newcommand{\qz}{{q_z}}

\newcommand{\fluxx}{{\varphi_x}}					 
\newcommand{\fluxy}{{\varphi_y}}
\newcommand{\tenseurxx}{{\varphi_{xx}}}
\newcommand{\tenseurxy}{{\varphi_{xy}}}

\newcommand{\e}{\varepsilon}
\newcommand{\ecarre}{{\varepsilon_2}}

\newcommand{\sfluxx}{s_\fluxx}					 
\newcommand{\sfluxy}{s_\fluxy}
\newcommand{\stenseurxx}{s_\tenseurxx}
\newcommand{\stenseurxy}{s_\tenseurxy}
\newcommand{\se}{s_{\e}}

\newcommand{\secarre}{s_\ecarre}

\newcommand{\sigmafluxx}{\sigma_\fluxx}					 
\newcommand{\sigmafluxy}{\sigma_\fluxy}
\newcommand{\sigmatenseurxx}{\sigma_{\tenseurxx}}
\newcommand{\sigmatenseurxy}{\sigma_{\tenseurxy}}
\newcommand{\sigmae}{\sigma_{\e}}

\newcommand{\sigmaecarre}{\sigma_{\ecarre}}

\newcommand{\eerho}{E_\e^\rho}					 
\newcommand{\eeqx}{E_\e^\qx}
\newcommand{\eeqy}{E_\e^\qy}
\newcommand{\eeqz}{E_\e^\qz}
\newcommand{\eecarrerho}{E_{\ecarre}^\rho}
\newcommand{\eecarreqx}{E_{\ecarre}^\qx}
\newcommand{\eecarreqy}{E_{\ecarre}^\qy}
\newcommand{\efluxxrho}{E_{\fluxx}^\rho}
\newcommand{\efluxxqx}{E_{\fluxx}^\qx}
\newcommand{\efluxxqy}{E_{\fluxx}^\qy}
\newcommand{\efluxyrho}{E_\fluxy^\rho}
\newcommand{\efluxyqx}{E_\fluxy^\qx}
\newcommand{\efluxyqy}{E_\fluxy^\qy}
\newcommand{\etenseurxxrho}{E_\tenseurxx^\rho}
\newcommand{\etenseurxxqx}{E_\tenseurxx^\qx}
\newcommand{\etenseurxxqy}{E_\tenseurxx^\qy}
\newcommand{\etenseurxyrho}{E_\tenseurxy^\rho}
\newcommand{\etenseurxyqx}{E_\tenseurxy^\qx}
\newcommand{\etenseurxyqy}{E_\tenseurxy^\qy}

\newcommand{\s}{\substack }					           
\newcommand{\dis}{\displaystyle}
\newcommand{\guill}[1]{\textquotedblleft{#1}\textquotedblright}
\newcommand{\pt}{\mbox{\begin{LARGE}$\cdot$\end{LARGE}}}

\newcommand{\qq}{\forall\,}
\newcommand{\ddqq}{D$d$Q$q$ }

\newcommand{\rg}{S(n,d)}
\usepackage{color}

\begin{document}

\title{Isotropy conditions for lattice Boltzmann schemes. \\ Application to D2Q9}
\thanks{This work has been financially supported by the French Ministry of Industry (DGCIS) and the Region Ile-de-France in the framework of the LaBS Project \cite{labs}.}
\author{Adeline Augier}
\address{Universit\'e Paris Sud, Laboratoire de math\'ematiques, UMR 8628, Orsay, F-91M05, France;\\ \email{adeline.augier@math.u-psud.fr,\ francois.dubois@math.u-psud.fr \&\ benjamin.graille@math.u-psud.fr} }
\author{Francois Dubois}
\sameaddress{1}
\secondaddress{
Conservatoire National des Arts et M\'etiers, Department of mathematics, Paris, France}
\author{Benjamin Graille}
\sameaddress{1}

\begin{abstract} 
In this paper, we recall the linear version of the lattice Boltzmann schemes in the framework proposed by d'Humi\`eres. According to the equivalent equations we introduce a definition for a scheme to be isotropic at some order. This definition is chosen such that the equivalent equations are preserved by orthogonal transformations of the frame. 
The property of isotropy can be read through a group operation and then implies a sequence of relations on relaxation times and equilibrium states that characterizes a lattice Boltzmann scheme.
We propose a method to select the parameters of the scheme according to the desired order of isotropy. Applying it to the D2Q9 scheme yields the classical constraints for the first and second orders and some non classical for the third and fourth orders. 
\end{abstract}


\maketitle
\section*{Introduction}
\label{sec:intro}

Numerous phenomena, in particular acoustic properties, are modelled by partial differential equations (PDEs) like Navier-Stokes equations. One of the main property of these PDEs is to take into account the isotropy of the environment. Namely, if the physical properties of the environment do not depend on the orientation, the associated PDEs are invariant by orthogonal transformations of the spatial frame.  
In this paper, we are interested by lattice Boltzmann schemes \cite{H94,LL00,D07} used by the industrialists because of their efficiency. As the directions of the lattice are privileged, they are not \textit{a priori} numerically isotropic. Because of the criss-cross pattern of the scheme \ddqq (where $d$ represents the space dimension and $q$ the number of discrete velocities), it is natural to look for parameters that give the same behaviour on each axes (\textit{i.e} by a rotation of angle $\pi/2$ in 2D for example) and by extent we want to obtain this result for each orthogonal transformation of the spatial frame. Our purpose is then to understand and control the lack of isotropy for these schemes. 

The first thing to do consists in defining the isotropy for lattice Boltzmann schemes. Actually, there is no uniform way to define isotropy: for example in \cite{CWSD92} the authors focus on the maintain of isotropy of the stress tensor required by fluid equations. An other point of view is described in \cite{LL00}: an isotropic fluid is simulated by eliminating as much as possible the effects due to the symmetry of the lattice and it means solving an eigenvalue problem. The problem of the space transformation is evoked in \cite{vdSE99,YGL05}: in \cite{YGL05} rotating frames are considered and in \cite{vdSE99} the authors introduce a definition of isotropic diffusion that constraints the collision operator to be invariant under all isometries of the Bravais lattice. In this paper, we define isotropy as the invariance by the orthogonal transformations of the frame and we express it on the equivalent equations. Namely, the equivalent equations are PDEs that come from a Taylor expansion of the lattice Boltzmann scheme \cite{D07} and a PDE is said isotropic if it is invariant by rewriting it in the \guill{new} frame. Finally, a \ddqq frame is isotropic at order $M$ if its set of equivalent equations of order $n$, $1\leq n\leq M$, is isotropic. 

In the first Section we quickly recall the lattice Boltzmann method and give some results on the equivalent equations.

In the second Section, we exploit the definition of isotropy for PDEs in an algebraic point of view. More precisely,  
we first define a new frame by orthogonal transformation. After redefining all the objects in the new frame, we identify both of the equivalent equations: the old one and the new one. The isotropy is guaranteed as soon as the various tensors in the PDEs appear as a fixed point of some group operation.
Then, we obtain a very useful set of relations on the coefficients of the equivalent equations.
 
In the last Section, a lattice Boltzmann scheme is said isotropic at order $M$ if its set of equivalent equations of order $n$, $1\leq n\leq M$, is isotropic in the sense explained in Section~\ref{sec:rotating}.
Since the coefficients of the equivalent equations depend on the parameters of the scheme (relaxation times and  equilibrium states), we are able to characterize,  for the related scheme,
the lack of isotropy on these parameters for each order.
More precisely, we are able to determine order after order the relations between the parameters of the scheme that have to be satisfied to have isotropic equivalent equations. We determine these relations for the D2Q9 scheme and specify all the sets of parameters that give isotropy until the fourth order. 

The lattice Boltzmann schemes considered in this paper are linear (for applications reasons) while the proposed definition of isotropy is more general and can be used as soon as the equivalent equations of the scheme are known.

\section{Lattice Boltzmann method and equivalent equations}
\label{sec:lbm}

In this section, we recall some notations related to the lattice Boltzmann method. Then, we introduce the equivalent equations as well as their link with the physic. This link is essential because it allows us to decide what parameters have physical values and what parameters have to be chosen in a allowed set of parameters. We will see below (in section~\ref{sec:isotropy}) what for a choice we have for this set of parameters in order to guaranty isotropy of the scheme.

\subsection{Notations and description of the lattice Boltzmann method}

We use the scheme proposed by d'Humi\`eres in \cite{H94}: we consider a regular lattice $\mathcal{L}$ with typical mesh size $\Delta x$. The time step $\Delta t$ is determined thanks to the velocity scale $\lambda$: $\Delta t=\Delta x/\lambda.$
For the \ddqq scheme, we note $\vectv=(v_j)_{1\leq j\leq q}$ the set of $q$ velocities and we assume that for each node $x$ of $\mathcal{L}$, and  each $v_j$ in $V$, the point $x+v_j\Delta t$ is also a node of the lattice $\mathcal{L}$. The aim of the \ddqq scheme is to precise the particle distribution $\vectf=\trans{(f_j(t,x))}_{1\leq j\leq q}$ (where $\trans{}$ represents the transpose) for $x\in\mathcal{L}$ and discrete values of time $t$ by solving a discretization of the Boltzmann's equation in two steps for every time: collision and transport. We finally have to propagate the so-called lattice Boltzmann scheme: $f_j(x+v_j\Delta t,t+\Delta t)\ = \ f_j^*(x,t),\ \forall\,x\in\mathcal{L}\,,\ 1\leq j\leq q$, where $f_j^*$ represents the $j$th component of the particle's distribution after the collision. In order to implement this scheme, as it is done in \cite{LL00}, we introduce the \textit{orthogonalized} moments $\vectm$ defined by $\vectm\ :=\ M\,\vectf$ where $M$ is a given and invertible matrix. Concerning the moments $\vectm$, we distinguish those who are preserved during the collision, called $\vectW\in\mathbb{R}^{N}$, of the complementary ones denoted by  $\vectY\in\mathbb{R}^{q-N}$.  
Since we are interested in linear applications, we assume that the collision operator is linear, so there exists a matrix $J$ such that $\vectm^*\ :=\ J\,\vectm$, with the following structure: 
$$\vectm^*\ =\ \left(\begin{array}{c}
         \vectW^*\\
	 \vectY^*
        \end{array}
\right)\ =\ J\,\vectm\ = \ \left(\begin{array}{cc}
   \Id_{N}&0\\
    S\,E&\Id_{q-N}-S
  \end{array}\right)\,\left(\begin{array}{c}
         \vectW\\
	 \vectY
        \end{array}\right),
$$
where $E\in\mathcal{M}_{q-N,N}(\mathbb{R})$ and $S$ is the diagonal matrix of the relaxation times $s_k=\Delta t/\tau_k$, $N+1\leq k\leq q$, with  $0< s_k<2$, $ N+1\leq k\leq q$ for stability conditions. 
We introduce a last classical notation: $\sigma_k\ =\ 1/s_k-1/2 \in\mathbb{R}^+,\ 1\leq k\leq q$.  
Finally a time step of linear lattice Boltzmann scheme scheme reduces to:
\begin{equation}
 \label{ddqq}
 \vectf(x+v_j\Delta t,t+\Delta t)\  =\ M^{-1}\,J\,M\,\vectf(x,t),\qquad \qq x\in\mathcal{L},\ 1\leq j\leq q. 
\end{equation}
Thus, the parameters of the lattice Boltzmann schemes are the coefficients of the matrices $S$ and $E$. This coefficients are usually chosen thanks to the experiment (by the physicists and the industrialists). The aim of this paper is to precise these coefficients according to the property of isotropy. 
\begin{rmrk}
There exist three types of lattice Boltzmann schemes according to the values of the relaxation times: either all of relaxation times are equal (BGK  for Bhatnagar - Gross - Krook scheme) or there is only two different relaxation times (TRT for Two Relaxation Times) or there is no conditions on these (MRT scheme for Multiple Relaxation Times). Since we investigate the behavior of the scheme in the sense of the isotropy, it is natural to consider the third case: it is the most general case even if it is not the most often used. Thus, there are $q-N$ independent unknowns parameters on $S$. 

Focus now on the matrix $E$ that is chosen full in order to give us as degree of freedom as possible: it gives $N(q-N)$ additional parameters. Actually, it is not an usual choice: the matrix $E$ is usually chosen sparse by physicists, but in our knowledge, there is no mathematical justification to this type of assumption (based on the physical \guill{good sense}). It is remarkable, that the isotropy condition for the \ddqq scheme can give some justification of this set of coefficients. 

Finally, both matrices $E$ and $S$ give $(N+1)\times(q-N)$ parameters to be precise in order to conserve physicals and isotropic properties. 
\end{rmrk}

For example, for the D2Q9 scheme, we consider the classical moments: the density $\rho$ and the momentum $\vectq=\trans{(\qx,\qy)}$ that are conserved during the collision and the energy $\e$, the square of the energy $\ecarre$, the flux $\trans{(\fluxx,\fluxy)}$, and the tensor $\trans{(\tenseurxx,\tenseurxy)}$. Then, we orthogonalize them like it is done in \cite{DL10}. For this choice of moments we define both equilibrium states and relaxation times matrices by:
$$\trans{E} \ := \
\left(\begin{array}{llllll}
\eerho\lambda^2\ & \eecarrerho\lambda^4\ & \efluxxrho\lambda^3 & \efluxyrho\lambda^3 & \etenseurxxrho\lambda^2 & \etenseurxyrho\lambda^2 \\
\eeqx\lambda\    & \eecarreqx\lambda^3\    & \efluxxqx\lambda^2  & \efluxyqx\lambda^2  & \etenseurxxqx\lambda    & \etenseurxyqx\lambda\\
\eeqy\lambda\    & \eecarreqy\lambda^3\    & \efluxxqy\lambda^2  & \efluxyqy\lambda^2  & \etenseurxxqy\lambda    & \etenseurxyqy\lambda\\
\end{array}\right)\ \mbox{and}\ S\ :=\ \mbox{Diag}(\se,\secarre,\sfluxx,\sfluxy,\stenseurxx,\stenseurxy).$$

\subsection{Equivalent equations}

We first introduce the notation $\pt^i_j$ for the element of $\pt$ that is on the intersection of the $i$th line and the $j$th column. By convention, a Latin letter is an index related on the moments (of size $N$) while Greek letter represents an index related on the space dimension (of size $d$). We then introduce the system of $N$ equivalent equations \cite{D07,D08} of order $M$ with which the linear lattice Boltzmann scheme is consistent at order $M$:
\begin{equation}
 \label{equivalent_equations}
\partial_t\vectW^i+ \sum_{n=1}^M\A{n}^i  \contr \nabla^n\vectW\ =\ 0,\ 1\leq i\leq N,
\end{equation}
 where $\A{n}^i,\ 1\leq n\leq M,\ 1\leq i\leq N$ are tensors of order $n$ that take into account the coefficients of $E$ and $S$. By convention, the maximal contraction operator $\contr$ is defined by: 
 $$\left(\A{n}\contr\nabla^n\vectW\right)^i\ :=\ \A{n}^i\contr\nabla^n\vectW\ :=\ \dis\sum_{\s{1\leq j\leq N\\1\leq \alpha_1,\cdots,\alpha_n\leq d}}\A{n}_j^{i,\alpha_1,\cdots,\alpha_n}\partial_{\alpha_1}\cdots\partial_{\alpha_n}\vectW^j,\qquad 1\leq i\leq N,\ 1\leq n\leq M. $$

These equivalent equations come from a formal calculus explained in \cite{DL10} and an algorithm that is easy to use is described in \cite{Dpreparation}.

\begin{rmrk}
The tensors $\A{n}$ belongs to the space $\mathbb{R}^{N^2n^d}$, $1\leq n\leq M$, but because of the Schwarz property, numerous of coefficients of $\A{n}^i$ are equal (or are proportional between them). In order to count how much of these coefficients are independent, we have to calculate the number of $d$-uplet $(\beta_1,\cdots,\beta_d)\in\mathbb{N}^d$ such that $\sum_{k=1}^d\beta_k=n$. This is a well-known combinatorics result: this number is the number of $d$ combinations with repetition of a set which cardinal number is $n$ and it reads: $\Gamma_d^n\ :=\ (d+n-1)!/(n!(d-1)!)$ and $\A{n}\in\rg$, where $\rg$ is the space of the symmetrical tensors of $\mathbb{R}^{N^2n^d}$. 
\end{rmrk}

The equivalent equations \eqref{equivalent_equations} are also very useful to make the link with the physics. For example, it is well-known (see for example \cite{D08}) that there exists a D2Q9 scheme (\textit{i.e.} a choice of moments $\vectY\in\mathbb{R}^6$) such that \eqref{equivalent_equations} is consistent with Navier-Stokes equations at the order 2. Then, there exists a relation between both viscosities (shear viscosity $\zeta$ and bulk viscosity $\mu$) and two relaxation times (see for example \cite{DL10}). Furthermore, equilibrium state of the energy given by $\eerho\rho+\eeqx\qx+\eeqy\qy+\eeqz\qz$ is such that $\eerho$ only depend on the sound velocity.  
 
Finally, our purpose is to precise as much as possible the $(N+1)(q-N)-3$ parameters of the lattice Boltzmann scheme in order to understand the lack of isotropy of these schemes: 
for example, in the D2Q9 scheme, parameters that are given thanks to the environment are $\eerho$ (it depends on the sound velocity), $\sigmae$ (it depends on the bulk viscosity $\zeta$) and $\sigmatenseurxx$ (it depends on the shear viscosity $\mu$) that means $21$ parameters.

\section{Algebra for rotating the equivalent equations}
\label{sec:rotating}

We first recall the definition of an isotropic PDE:
\begin{dfntn}
\label{def:isotropyPDE}
 A PDE  
\begin{equation}
 \label{PDE}
\partial_t\vectW^i+ \sum_{n=1}^M\A{n}^i\contr\nabla^n\vectW\ =\ 0,\ 1\leq i\leq N,
\end{equation}
is said isotropic if it is invariant by orthogonal transformation of the frame. 
\end{dfntn}

We introduce the matrix $r\in\mathcal{O}_d(\mathbb{R})$ that represents an orthogonal transformation of the frame $(O,\vecte)$. We then use the notation $\widetilde{\pt}$ for every element in the \guill{new} frame: in this way, if the vectors $\vectx\in\mathbb{R}^d$ and the moments $\vectW\in\mathbb{R}^N$ are known in the frame $(O,\vecte)$, we have to define the news moments $\vectWt$, the new vector $\vectxt$ (such that $\vectx:=r\,\vectxt$), the new basis $\vectet$, the news tensors $\At{n}^i$, $1\leq n\leq M$, $1\leq i\leq N$ and the new partial derivatives $\widetilde{\partial_\alpha}$, $1\leq\alpha\leq d$. Finally, let be $R(r)$ the matrix in $\mathcal{O}_{N}(\mathbb{R})$ such that $\vectW\ :=\ R(r)\,\vectWt$. More precisely, $R: \mathcal{O}_{d}(\mathbb{R})\rightarrow\mathcal{O}_{N}(\mathbb{R})$ defines a morphism that gives a group representation. In the present contribution, mass and momentum are conserved thus $\vectW=\trans{(\rho,\vectq)}$, $N=d+1$ and for every $r\in\mathcal{O}_d(\mathbb{R})$, we get $R(r)=\mbox{Diag}(1,r)\in\mathcal{O}_N(\mathbb{R}$). 

The aim of this section consists in giving the relation between objects that are defined before and after the transformation of the frame. First, we give a relation between both partial derivatives and we deduce a relation on both derivatives of the moments. Second, we introduce a group operation \textit{via} an homomorphism that will show itself very useful for the definition of the isotropy. 

\begin{lmm}
For every $1\leq \alpha\leq d$ and $1\leq i\leq N$, the derivatives $\partial_\alpha\vectW^i$ at the $i$th component of the moment $\vectW$ reads:
\begin{equation}
\label{eq:derivatives_moments}
\partial_\alpha\vectW^i\ =\ \partial_\alpha(R(r)\vectWt)^i\ =\ \sum_{j,\beta}R_j^i(r)(r^{-1})_\alpha^\beta\widetilde{\partial_\beta}\vectWt^j.
\end{equation}
\end{lmm}
\begin{proof}
As $\vectx={r}\,\vectxt$, for each $1\leq \alpha\leq d$, we have the following relation on the partial derivatives $\partial_\alpha$ along the $\alpha$th direction:
\begin{equation*}
\label{eq:derivatives}
 \partial_\alpha=\sum_{1\leq \beta\leq d}\dfrac {\partial\vectxt^\beta}{\partial\vectx^\alpha}\cdot\widetilde{\partial_\beta}=\sum_{1\leq \beta\leq d}\left(r^{-1}\right)^\beta_\alpha\widetilde{\partial_\beta}.
\end{equation*}
Then, it is easy to characterize the partial derivatives of the moments and we obtain the desired result \eqref{eq:derivatives_moments}.
\end{proof}
 We also define a group operation through the following Proposition.
\begin{dfntnprpstn}
 Let $n$ be in $\mathbb{N}^*$ and the map $\Phi_n\ :\ \mathcal{O}_d(\mathbb{R})\times\rg\rightarrow\rg$ defined for $1\leq i,j\leq N$, $1\leq\alpha_1,\cdots,\alpha_n\leq d$, $1\leq n\leq M$, by the relation:
\begin{equation}
\label{def:Phi_n}
\left(\Phi_n(r)(\A{n})\right)^{i,\alpha_1,\cdots,\alpha_n}_j\ :=\ \sum_{\s{1\leq \beta_1,\cdots,\beta_n\leq d\\1\leq k,l\leq N}}\left(R(r)\right)^i_l\A{n}^{l,\beta_1,\cdots,\beta_n}_kr^{\alpha_1}_{\beta_1}\cdots r^{\alpha_n}_{\beta_n} (R(r)^{-1})^k_j.
\end{equation}
It is an homomorphism that characterizes a group operation. 
\end{dfntnprpstn}

\begin{proof}
 Let $\Id$ be the identity matrix of $\mathcal{O}_d(\mathbb{R})$, then for every $\A{n}\in\rg$ we have $\Phi_n(\Id)(\A{n})=\A{n}$. It remains to prove the following property: $\qq r,r'\in\mathcal{O}_d(\mathbb{R})$, $\qq \A{n}\in\rg$, $\Phi_n(rr')(\A{n})=\Phi_n(r)\left(\Phi_n(r')(\A{n})\right)$. First, by using the definition of the homomorphism $\Phi_n$ given in~\eqref{def:Phi_n}, we get for $1\leq i,j\leq N$ and $1\leq \alpha_1,\cdots,\alpha_n\leq d$:
\begin{multline}
\label{dem_phi_interm}
 \Phi_n(r)\left(\Phi_n(r')(\A{n})\right)^{i,\alpha_1,\cdots,\alpha_n}_j\ \\
=\ \sum_{\s{1\leq \beta_1,\cdots,\beta_n,\gamma_1,\cdots,\gamma_n\leq d\\1\leq k,l,m,o\leq N}}R(r)^i_lR(r')^l_m\,\A{n}^{m,\gamma_1,\cdots,\gamma_n}_or^{\alpha_1}_{\beta_1}(r')_{\gamma_1}^{\beta_1}\cdots r^{\alpha_n}_{\beta_n}(r')_{\gamma_n}^{\beta_n}(R(r')^{-1})^o_k(R(r)^{-1})^k_j.
\end{multline}
Since $R(rr')=R(r)R(r')$, and so $R(rr')^{-1}=R(r')^{-1}R(r)^{-1}$, \eqref{dem_phi_interm} gives: $\Phi_n(r)\left(\Phi_n(r')(\A{n})\right)\ = \ \Phi_n(rr')(\A{n})$
and the homomorphism $\phi_n$ truly defines a group operation.
\end{proof}

\begin{prprt}
 Let $r$ be an orthogonal transformation of $\mathcal{O}_d(\mathbb{R})$, then the system \eqref{PDE} reads in the new frame 
\begin{equation}
 \label{PDEr}
\partial_t \vectWt^i+ \sum_{n=1}^M\At{n}^i\contr\nablat^n\vectWt=0,\quad 1\leq i\leq N,
\end{equation}
where $\A{n}\ =\ \Phi_n(r)(\At{n})$, $1\leq n\leq M$, or equivalently $\At{n}\ =\ \Phi_n(r^{-1})(\A{n})$, $1\leq n\leq M$.
\end{prprt}
\begin{proof}
Let us write \eqref{PDE} after the orthogonal change of frame $r$. Because of the equality $\vectW=R(r)\,\vectWt$, $\partial_t\vectW$ becomes $\partial_tR(r)\,\vectWt$. Using the property~\eqref{eq:derivatives_moments}, we obtain a relation that gives $\A{n}^i\contr\nabla^n\vectW$ in function of the partial derivatives of the moments in the new frame, more precisely for $1\leq i\leq N$ we get:
$$
\A{n}^i\contr\nabla^n\vectW\ =\ \sum_{\s{1\leq \beta_1,\cdots,\beta_n,\gamma_1,\cdots,\gamma_n\leq d\\1\leq k,l,m,o\leq N}}\A{n}^{i,\alpha_1,\cdots,\alpha_n}_j\left(r^{-1}\right)_{\alpha_1}^{\beta_1}\cdots \left(r^{-1}\right)_{\alpha_n}^{\beta_n}R(r)^j_k\widetilde{\beta_1}\cdots\widetilde{\beta_n}\vectWt^k.
$$
Thus, by identifying both \eqref{PDE} and $R(r)\times\eqref{PDEr}$, we characterize $\At{n}$, $1\leq n\leq M$, thanks 
$$\At{n}^i\contr\nablat^n\vectWt= \sum_{1\leq j\leq N}\left(R(r)^{-1}\right)^i_j\,\left(\A{n}^j\contr\nabla^n\vectW\right),\ \qq \vectW\in\mathbb{R}^{N},\ 1\leq n\leq M,\ 1\leq i\leq N,$$ and that gives the desired result.   
\end{proof}

Finally, we obtain a characterization for a PDE to be isotropic:
\begin{crllr}
 Let \eqref{PDE} be a PDE of order $M$, it is isotropic if and only if $\A{n}\ =\ \At{n}$, $1\leq n\leq M$, $\qq r\in\mathcal{O}_d(\mathbb{R})$ and it means that the tensor $\A{n}$ is a fixed point of $\Phi_n$, $1\leq n\leq M$, that is $\Phi_n(r)(\A{n})=\A{n}$, $\qq r\in\mathcal{O}_d(\mathbb{R})$, $1\leq n\leq M$. 
\end{crllr}

\section{Isotropy conditions for lattice Boltzmann schemes}
\label{sec:isotropy}

In this section, 
we first give a definition of isotropy for lattice Boltzmann schemes and then we give 
an application on the D2Q9 scheme.  
\begin{dfntn}
\label{def:isotropy}
A lattice Boltzmann scheme is said isotropic at the order $M$ if the system of equivalent equations \eqref{equivalent_equations} at the order $M$ is isotropic. 

Furthermore, we denote by $L_N(r)\ :=\ \sum_{1\leq n\leq N}(\Phi_n(r)(\A{n})-\A{n}){\triangle t}^n$ the lack of isotropy at order $N$ of the orthogonal transformation $r$ of the equivalent equations \eqref{equivalent_equations}.
\end{dfntn}

Then, we use this definition on the D2Q9 scheme in order to obtain all the sets of parameters that improve this scheme (in the sense of the isotropy) until the fourth order. 
In this case, the orthogonal transformation $r$ is the well-known rotation matrix, $d$ is equal to $2$, $N=3$ and $M=4$. 

\begin{prpstn}
\label{prop}
 Let $L_4(r)$ be the lack of isotropy for the D2Q9 scheme at fourth order for the orthogonal transformation $r$, then we get:
\begin{itemize}
 \item $L_4(r)=0({\Delta t}^2),\qq r\in\mathcal{O}_2(\mathbb{R})$ iff $\eeqx=\eeqy=\etenseurxxrho=\etenseurxxqx=\etenseurxxqy=\etenseurxyrho=\etenseurxyqx=\etenseurxyqy=0$.
 \item $L_4(r)=0({\Delta t}^3),\qq r\in\mathcal{O}_2(\mathbb{R})$ 
iff $\efluxxrho=\efluxxqy=\efluxyrho=\efluxyqx=0$ and $\efluxxqx=\efluxyqy$.
 \item $L_4(r)=0({\Delta t}^4),\qq r\in\mathcal{O}_2(\mathbb{R})$ 
iff (\/$\eecarreqx=\eecarreqy=0$, $\sigmatenseurxx=\sigmatenseurxy$, $\efluxxqx=-1$\/) and (either $2\eecarrerho+4+3\eerho=0$ or $\sigmafluxx=\sigmafluxy=1/(12\sigmatenseurxx)$\/).
 \item $L_4(r)=0({\Delta t}^5),\qq r\in\mathcal{O}_2(\mathbb{R})$ 
iff (\/$2\eecarrerho+4+3\eerho=0$, $\sigmae=\sigmatenseurxx$, $\sigmafluxx=\sigmafluxy=1/(6\sigmatenseurxx)$\/) and (either $2+3\eerho=0$ or $\sigmaecarre=\sigmatenseurxx$\/). Thanks to \cite{DL10} this last condition gives the following relation on both viscosities: $\zeta=6\mu$.
\end{itemize}
\end{prpstn}
\begin{rmrk}
Some of the properties for an isotropic D2Q9 scheme that are obtained here are well-known and naturally used in this scheme. In fact, the parameters that are taken null in order to obtain isotropy at first and second orders are also taken null for isotropic reasons considering the kinetic solution at equilibrium in the continuous environment.
However, the results on the third and fourth orders are more surprising even if some usual D2Q9 schemes can be seen as particular cases of those we propose.
\end{rmrk}
 
The proof of Proposition~\ref{prop} is obtained thanks to a rigorous implementation of isotropy conditions with a software of formal calculus and \cite{ADGpreparation} give details of it.

\bibliographystyle{plain}
\nocite{}

\bibliography{bibliographie_lbm}
\end{document}